\documentclass[12pt]{amsart}

\usepackage{amssymb, latexsym, mathtools}
\input xy
\xyoption{all}

\def\mbR{\mathbb{R}}

\theoremstyle{plain}
\newtheorem{theorem}{Theorem}[section]

\newtheorem{proposition}[theorem]{Proposition}
\newtheorem{lemma}[theorem]{Lemma}

\theoremstyle{definition}
\newtheorem{definition}[theorem]{Definition}
\newtheorem{remark}[theorem]{Remark}
\newtheorem{example}[theorem]{Example}
\newtheorem{problem}[theorem]{Problem}

\numberwithin{equation}{section}

\usepackage
{hyperref}

\usepackage{braket}

\usepackage{tikz}
\usepackage{tkz-euclide}
\usepackage{tikz-3dplot}
\usepackage{graphicx}

\begin{document}

\title[Barycentric coordinates]{From affine to barycentric coordinates in polytopes}

\author[Romanowska]{A.B. Romanowska$^1$}
\author[Smith]{J.D.H. Smith$^2$}
\author[Zamojska-Dzienio]{A. Zamojska-Dzienio$^3$}
\address{$^{1,3}$ Faculty of Mathematics and Information Science\\
Warsaw University of Technology\\
00-661 Warsaw, Poland}

\address{$^2$ Department of Mathematics\\
Iowa State University\\
Ames, Iowa, 50011, USA}

\email{$^1$anna.romanowska@pw.edu.pl\phantom{,}}
\email{$^2$jdhsmith@iastate.edu\phantom{,}}
\email{$^3$anna.zamojska@pw.edu.pl}

\urladdr{$^2$
\protect
{
\href{https://jdhsmith.math.iastate.edu/}
{https://jdhsmith.math.iastate.edu/}
}
}

\maketitle

\begin{abstract}
Each point of a simplex is expressed as a unique convex combination of the vertices. The coefficients in the combination are the barycentric coordinates of the point. For each point in a general convex polytope, there may be multiple representations, so its barycentric coordinates are not necessarily unique. There are various schemes to fix particular barycentric coordinates: Gibbs, Wachspress, cartographic, etc.  In this paper, a method for producing sparse barycentric coordinates in polytopes will be discussed. It uses a purely algebraic treatment of affine spaces and convex sets, with barycentric algebras. The method is based on a certain decomposition of each finite-dimensional convex polytope into a union of simplices of the same dimension.
\end{abstract}

\section{Introduction}

The first part of the paper is a brief survey of the main properties of the abstract algebras that respectively define affine spaces and convex sets, the latter generating a variety of non-associative algebras called \emph{barycentric algebras}. Section~\ref{S:baralg} provides the basic definitions and properties of such algebras. Sections~\ref{S:walls} and~\ref{S:simpl} deal with two important concepts concerning (convex) polytopes \textemdash convex sets generated as algebras by a finite number of elements. Section~\ref{S:walls} discusses walls (faces) of polytopes, and their role in the structure theory of barycentric algebras, while Section~\ref{S:simpl} provides basic properties of simplices considered as barycentric algebras, and in particular their role as free barycentric algebras. An important fact used later is that a nonempty set of $(n+1)$ affinely independent elements of a real affine space generates both the free affine space $\mathbb R^n$ and the corresponding free barycentric algebra, the $n$-dimensional simplex $\Delta_n$.

The new results of the paper are presented in Sections~\ref{S:barcordpolytope} and \ref{S:barcordDec}, which are devoted to our main problem:
\begin{quote}
Provide a barycentric coordinate system for points $\mathbf x$ of a convex polytope $P$.
\end{quote}
The coordinates express each point as a convex combination of certain vertices of the polytope. There are different approaches to the problem, each taken in various contexts. Section~\ref{S:barcordpolytope} and the references summarize many known results (see \cite{FloaterGBCA}, for example). Here, we mention three specific barycentric coordinate systems:
\begin{itemize}
\item
the most general \emph{Gibbs coordinates} based on entropy maximization, which actually work for all barycentric algebras \cite{HormSuku}, \cite[Ch.~IX]{Modes};
\item
the well-known \emph{Wachspress coordinates} for polytopes \cite{Wachspress2,WarSchHirDes}; and
\item
our recently formulated \emph{chordal} and \emph{cartographic coordinates} for polygons \cite{RSZ23}.
\end{itemize}
After formulating the problem, and giving some technical details on the relationship between affine and barycentric coordinates, we describe our approach to its solution. Chordal coordinates are based on chordal decompositions of a polygon into a union of triangles ($2$-dimensional simplices), and analysis of parsing trees of such decompositions. The method produces sparse coordinates that involve a minimal number of vertices of the polygon, and has interesting combinatorial properties. However, it does not carry over to polytopes of higher dimension.

In Section~\ref{S:barcordDec}, we describe our new method based on a special decomposition of an arbitrary polytope $P$ into a union of simplices of the same dimension as $P$. Each point $\mathbf x$ of $P$ is contained in a simplex of the decomposition, and has a unique expression as a convex combination of the vertices of this simplex. The barycentric coordinates of $\mathbf x$ in $P$ are the coordinates of this combination, with zero coordinates at the remaining vertices of $P$. In this system, the barycentric coordinates of a point $\mathbf x$ are referred to a simplex $S$ of the decomposition precisely when the affine coordinates of $\mathbf x$ with respect to $S$ are all non-negative. We provide a method of numbering the simplices of a decomposition of $P$, and show how to calculate the barycentric coordinates of points of $P$.

\section{Barycentric algebras}\label{S:baralg}

As background, this section records basic facts concerning affine spaces, convex sets and barycentric algebras. Further details may be found in the monographs~\cite{AB83} (for the geometric approach) and~\cite{Modes} (for the algebraic approach). The survey~\cite{R18} and its bibliography cover further sources.

\subsection{Affine subspaces and convex subsets of real spaces}\label{S:affconv}

Write $I=[0,1]$ for the closed real unit interval, and $I^{\circ}$ for its interior $]0,1[$.\footnote{The notation $(0,1)$ is reserved for the ordered pair consisting of $0$ and $1$.} For elements $x,y$ of a real vector space $V$, we set
\begin{equation}\label{E:operp}
xy\, \underline{p} \coloneqq x(1-p) + yp
\end{equation}
for $p\in\mathbb R$. In particular, note $xx\,\underline p=x=xy\,\underline 0$ and $xy\,\underline 1=y$. Sometimes, we write $p' \coloneqq 1-p$ for $p \in I^{\circ}$. Then~\eqref{E:operp} may be written as
$
xy\underline{p} = xp'+ yp.
$

If two points $x,y$ of the vector space $V$ are distinct, they lie on the line
$
L_{x,y} = \set{xy\, \underline{p}| p \in \mathbb{R}}.
$
A subset $A$ of $V$ is an \emph{affine subspace} of $V$ if, together with any two distinct points $x$ and $y$ in $A$, we have $L_{x,y}\subseteq A$. The line segment $I_{x,y}$ joining distinct points $x, y$ of $V$ is described as
$
I_{x,y} = \set{xy\, \underline{p}|p \in I}.
$
A subset $C$ of $V$ is a \emph{convex subset} of $V$ if, together with any two distinct points $x$ and $y$ in $C$, we have $I_{x,y}\subseteq C$. Note that these definitions admit the empty set as an affine space and as a convex set.

\subsection{Affine spaces}\label{S:affspaces}

The observations of \S\ref{S:affconv} provide an interpretation of affine spaces as \emph{algebras},\footnote{This is the usage of general or universal algebra, not to be confused with the linear algebraic notion of a compatible ring structure on a vector space.} sets with a set of operations defined on them. Thus, affine spaces become algebras $(A,\underline{\mathbb{R}})$, where $\underline{\mathbb{R}} \coloneqq \{\underline{p} \mid p \in \mathbb{R}\}$.
Affine spaces defined in this way have all the properties required for affine spaces defined in the traditional way. (See e.g. \cite[Sec.~1.6]{Modes}.)
In particular, the affine subspaces of an affine space $A$ are precisely the subalgebras of $(A,\underline{\mathbb{R}})$ (subsets of $A$ closed under the operations $\underline{p}$), affine mappings $h\colon A \rightarrow A'$ are algebra homomorphisms $h\colon (A,\underline{\mathbb{R}}) \rightarrow (A',\underline{\mathbb{R}})$ (mappings preserving the operations $\underline{p}$); and the direct product of a family of affine spaces is a direct product of algebras with operations $\underline{p}$ defined componentwise. For a natural number $k$, a $k$-dimensional affine space is affinely isomorphic to the affine space $\mathbb{R}^k$. The empty set is considered as an affine space of dimension $-1$.

Elements of the affine hull of a set $\{x_0, x_1,\dots,x_n\}$ of affinely independent elements are described as \emph{affine combinations}
\begin{equation}\label{E:affcomb}
x = \sum_{i=0}^{n} x_{i} p_i\ \text{with} \ \sum_{i=0}^{n}p_{i} = 1,
\end{equation}
and with \emph{affine coordinates} $p_i \in \mathbb{R}$. The affine combinations form an affine space that is isomorphic to $\mathbb{R}^n$.
They may be obtained by composing binary affine combinations, the basic binary operations $\underline{p}$~\cite[Lemma~1.6.4.1, Cor.~6.3.5]{Modes}.
The class $\mathbf{A}$ of (algebras isomorphic to) affine spaces considered as algebras $(A,\underline{\mathbb{R}})$ is a variety. (See e.g. \cite[Ch.~6]{Modes}.) This means that $\mathbf{A}$ is closed under the formation of subalgebras, direct products and homomorphic images. As each variety of similar algebras may be equivalently described  as an equationally defined class of algebras, we can see that the variety of affine spaces is equationally defined. (See~\cite[Th.~6.3.4]{Modes} for the identities defining this variety.)

\begin{remark}
Note that each of the nontrivial (i.e., $p\notin\set{0,1}$) operations $\underline{p}$ of an affine space is in fact a \emph{quasigroup} operation: For given $a,b$, the equations $ax \underline{p} = b$ and $ya \underline{p} = b$ have respective unique solutions $x$ and $y$. Thus, affine spaces may be considered as multi-quasigroups.
\end{remark}

\subsection{Convex sets}\label{S:convsets}

In parallel with the algebraic treatment of affine spaces presented in \S\ref{S:affspaces}, the observations of \S\ref{S:affconv} allow us to treat convex sets as algebras $(C,\underline{I}^{\circ})$, where $\underline{I}^{\circ}\coloneqq\{\underline{p} \mid p \in I^{\circ}\}$.
In this context, we also say that a convex set is a \emph{subreduct} of an affine space. The smallest such affine space is called the \emph{affine hull} of the convex set. If this affine hull has finite dimension $k$, it is construed as the topological space $\mathbb{R}^k$ with the usual Euclidean topology. For the general case (which is not needed for the work of this paper), see \cite[Ex.~5.8K]{Modes}. A convex set $C$ is $k$-\emph{dimensional} for a natural number $k$ if its affine hull is $\mathbb{R}^k$. The empty convex set has dimension $-1$.

As algebras, \emph{convex polytopes} (or briefly \emph{polytopes}) are defined as finitely generated convex sets. The minimal set of generators of a polytope $C$ is the set of its vertices (i.e., its extreme points). In geometric terminology, the convex set generated by a set $X$ is its \emph{convex hull}. If a $k$-dimensional polytope $C$ has $n$ vertices, then $n$ is at least $k+1$. A polytope $C$ is closed when considered as a subset of its affine hull.

The class $\mathbf{C}$ of (algebras isomorphic to) convex sets considered as algebras $(C,\underline{I}^{\circ})$ is closed under the formation of subalgebras and direct products, but is not closed under homomorphic images. Thus, it is not a variety.

\begin{definition}\label{D:BarycAlg}
A \emph{real barycentric algebra} (or briefly a \emph{barycentric algebra}) is defined as a homomorphic image of a convex set.
\end{definition}

\begin{theorem}\cite[\S5.8]{Modes}
The class $\mathbf{B}$ of barycentric algebras is the variety axiomatized by
the identities
\begin{equation}\label{E:idem}
x x\, \underline{p} = x
\end{equation}
of \emph{idempotence} for $p\in I^\circ$, the identities
\begin{equation}\label{E:skewcom}
x y\, \underline{p} = yx\, \underline{p}'
\end{equation}
of \emph{skew-commutativity} for $p\in I^\circ$, and the identities
\begin{equation}\label{E:skewass}
[x y\, \underline{p}]\ z\, \underline{q} = x\,\, [y z\,
\underline{q/(p\circ q)}]\,\, \underline{p\circ q}
\end{equation}
of \emph{skew-associativity} for $p,q\in I^\circ$, where $p \circ q:= p+q-pq$.

The class $\mathbf{C}$ of convex sets is defined, within the variety $\mathbf{B}$, by
the quasi-identities
\begin{equation}\label{E:cans}
(xy\, \underline{p} = xz\, \underline{p}) \Rightarrow (y = z)\quad\mbox{and}\quad (xy\, \underline{p} = zy\, \underline{p}) \Rightarrow (x = z)
\end{equation}
of \emph{cancellativity} for $p \in I^{\circ}$.
\end{theorem}

Barycentric algebras also satisfy the identities
\begin{equation}\label{E:entr}
[xy\underline{p}]\,[zt\underline{p}]\,\underline{q} =
[xz\underline{q}]\,[yt\underline{q}]\,\underline{p}
\end{equation}
of \emph{entropicity} for $p, q \in I^{\circ}$.
Moreover, idempotence and entropicity imply the identities
\begin{equation}\label{E:distr}
[xy\underline{p}]\,z\, \underline{q} = [xz\underline{q}]\,[yz\underline{q}]\, \underline{p}\quad\mbox{and}\quad
x\,[yz\underline{p}]\, \underline{q} = [xy\underline{q}]\,[xz\underline{q}]\, \underline{p}
\end{equation}
of \emph{distributivity} for $p, q \in I^{\circ}$.
The entropic and distributive laws are satisfied by all affine spaces.

\begin{remark}
For certain purposes, the closed unit interval $\underline{I} \coloneqq \{\underline{p} \mid p \in I = [0,1]\}$ may be used as the set of basic operations when defining convex sets, instead of the open unit interval $\underline{I}^{\circ}$ (see \cite{RSO07}, for example).
However, our current choice of the open unit interval allows one to develop an algebraic structure theory of barycentric algebras which would not be possible using the closed unit interval. See Remark~\ref{R:RegulrId} for further elaboration of this important point.
\end{remark}

In what follows, we combine the algebraic and geometric approaches to affine spaces and convex sets.

\subsection{The structure of barycentric algebras}\label{S:structure}

By the results of the preceding section, the variety $\mathbf{B}$ of barycentric algebras certainly contains all convex sets within its subclass $\mathbf C$ comprising those barycentric algebras which satisfy the cancellation property \eqref{E:cans}. On the other hand, extreme examples of barycentric algebras which violate \eqref{E:cans} are obtained when all the operations of $\underline{I}^{\circ}$ are equal, i.e. $x \circ y \coloneqq xy\underline{p} = xy\underline{q}$ for all $p, q \in I^{\circ}$. In this case, skew-commutativity and skew-associativity of the operations from $\underline{I}^{\circ}$, along with their idempotence, imply that the operation $\circ$ is idempotent, commutative and associative. Idempotent, commutative semigroups are \emph{semilattices}. Semilattices considered as algebras $(A,\underline{I}^{\circ})$ are called (\emph{iterated}) \emph{semilattices}.

A semilattice $(S,\circ)$ may be construed as an ordered set $(S,\leq_\cdot)$, where for any $a, b \in S$, $a \leq_{\cdot} b$ iff $a \circ b = a$, and such that any two elements of $S$ have a greatest lower bound, called their \emph{meet}, and denoted by $a \cdot b$. Then $a \cdot b = a \circ b$.  Dually, $(S,\circ)$ can be defined as an ordered set $(S,\leq_{+})$, where for any $a, b \in S$, $a \leq_{+} b$ iff $a \circ b = b$, and such that any two elements of $S$ have a least upper bound, called their \emph{join}, and denoted by $a + b$. Then $a + b = a \circ b$. A \emph{lattice} is an algebra $(L,+,\cdot)$ with two semilattice operations, where both semilattice orderings coincide.

\begin{proposition}\cite{N70}, \cite[\S7.6]{Modes}
Iterated semilattices form the only non-trivial proper subvariety $\mathbf{Sl}$ of $\mathbf{B}$. It is defined relative to $\mathbf{B}$ by the identities
$$
xy \underline{p} = xy \underline{q}
$$
for all $p, q \in I^{\circ}$.
\end{proposition}

More complex barycentric algebras are constructed from convex sets and semilattices, as certain disjoint sums of convex sets. The starting point is the following basic lemma. (For details see~\cite[Chs.~3,~7]{Modes}.)

\begin{lemma}\label{L:slqt}
Each barycentric algebra has a homomorphism onto an iterated semilattice.
\end{lemma}

The largest semilattice quotient $S$ of a barycentric algebra $A$ is called its \emph{semilattice replica}. If $\rho\colon A \rightarrow S$ is the \emph{semilattice replication homomorphism} of $A$ onto its semilattice replica $S$, then for each $s \in S$, the set $\rho^{-1}\set s$ is a subalgebra of $A$, and is a convex set. In such a context, we say that the barycentric algebra $A$ is a \emph{semilattice sum} of the $S$-ordered system $\set{\rho^{-1}\set s|s\in S}$ of convex sets. Each of these convex sets $\rho^{-1}\set s$ is an open subset of its affine hull \cite[\S~7.5]{Modes}. We obtain the following.

\begin{theorem}\label{T:strthm}\cite[Cor.~7.5.9]{Modes}
Each barycentric algebra is a semilattice sum of open convex sets.
\end{theorem}

The semilattice underlying Theorem~\ref{T:strthm} is the semilattice replica of the barycentric algebra. Theorem~\ref{T:strthm} is the basis of various structure theorems for barycentric algebras, provided e.g. in \cite[\S7.5]{Modes}, \cite{RS90, RS93}. Some other characterizations of convex sets and barycentric algebras are also presented there.

\begin{remark}\label{R:RegulrId}
Lemma~\ref{L:slqt} shows that the identities true in barycentric algebras must also be satisfied by all (iterated) semilattices. It is known that the semilattice identities are all regular, i.e. each has the same sets of variables on both sides. (See~\cite[Ex~1.5.4]{Modes}.) Thus, the identities true in barycentric algebras are also regular. On the other hand, the projections $\underline{0}$ and $\underline{1}$ satisfy irregular identities. As a consequence, barycentric algebras defined with the closed unit interval $I$ only have trivial semilattice replicas, and structure theorems based on Theorem~\ref{T:strthm} are not available in this situation.
\end{remark}

 \begin{example}\label{Ex:int}
Consider the real closed unit interval $I$ as a barycentric algebra, and its semilattice replication homomorphism $\rho\colon I \rightarrow S$. The semilattice replica $S$ has three elements $a > b < c$ with $\rho^{-1}\set a = \{0\}$, $\rho^{-1}\set c = \{1\}$ and $\rho^{-1}\set b = I^{\circ}$. The barycentric algebra $I$ is the (meet) semilattice sum of two one-element convex sets and the open unit interval $I^\circ$ over the semilattice $S$.
\end{example}

\begin{example}\label{Ex:extreal}
Consider the barycentric algebra $\mathbb R^\infty$ of real numbers extended by $\infty$. Here, $\mathbb R$ appears as a subalgebra with its convex set structure, while $x\infty\underline p = \infty$ for all $x\in\mathbb R^\infty$ and $p\in I^\circ$. The semilattice replica $S$ of $\mathbb R^\infty$ is the two-element (join) semilattice $\mathbf{2}=\set{a < b}$. The semilattice replica homomorphism $\rho\colon \mathbb R^\infty \rightarrow \mathbf{2}$ assigns $a$ to all elements of $\mathbb{R}$ and $b$ to $\infty$. The algebra
$\mathbb R^\infty$ is the semilattice sum of two convex sets: $\mathbb{R}$ and $\{\infty\}$.
\end{example}

We conclude this section with a structure theorem for barycentric algebras which is not directly related to semilattice sums of convex sets.

\begin{theorem}\label{T:structurethm}\cite[Th.~7.5.11]{Modes}
Each barycentric algebra $A$ embeds into a product of copies of the barycentric algebra $\mathbb{R}^{\infty}$.
\end{theorem}

\section{Walls}\label{S:walls}

Suppose that $(P,\underline I^\circ)$ is a barycentric algebra. A subset $A$ of $P$ is a \emph{wall} if
$$
\forall\ a,b\in P\,,\  \forall\,\ r \in I^\circ\,,\
a b\underline{r} \in A
\
\Leftrightarrow
\
a\in A
\mbox{ and }
b\in A\,.
$$
(The concept of a wall readily extends to general abstract algebras \cite[Def'n.~1.1.9(b)]{Modes}.)
If $P$ is a polytope, the algebraic concept of a wall matches the geometric concept of a face (see e.g. \cite[\S5, \S7]{AB83}). The trivial or empty subset of $P$ is the unique $(-1)$-dimensional face. The $0$-dimensional faces of an $n$-dimensional polytope $P$ are its vertices, and the $1$-dimensional faces are the edges. The maximal (i.e., $(n-1)$-dimensional) faces are called the \emph{facets} of $P$. The faces of a polytope are again polytopes, and each face of $P$ considered as a barycentric algebra is generated by the vertices of $P$ that it contains. The \emph{boundary} of $P$ in the affine hull of $P$ is the difference between $P$ and its (relative) interior, and is the union of the proper faces of $P$ (cf. \cite{LeeRecent}.)

The faces of a polytope $P$ form a lattice under inclusion. The meet of two faces $E$ and $F$ is their intersection (which may well be empty), while the join is the smallest face that contains $E$ and $F$. The empty face $\O$ is the smallest element of this lattice, and the entire polytope $P$, the unique improper face, is the largest element.

For a barycentric algebra $A$ and element $a \in A$, let $[a]$ be the \emph{principal wall} generated by $a$, i.e., the smallest wall containing $a$. The join of two principal walls is again a principal wall \cite[Lemma~3.3.9]{Modes}.
The following proposition is a special case of a more general result \cite[Th.~3.3.11]{Modes}.

\begin{proposition}\label{P:decomposition}
Principal walls of a barycentric algebra $A$ form a semilattice. The semilattice replica $S$ of $A$ is its semilattice of principal walls, and
$
\rho\colon A\to S;a\mapsto[a]
$
is the replication homomorphism.
\end{proposition}

Note that the nonempty faces of a polytope are the principal walls \cite[\S~5]{AB83}. As a principal wall, each face is generated by any of its interior points. In particular, any polytope $P$ is a semilattice sum of open convex sets (the interiors of the nonempty faces) over the (join) semilattice of its non-empty faces.

\section{Simplices}\label{S:simpl}

Simplices are polytopes, and each face of a simplex is itself a simplex. In geometrical terms, an $n$-dimensional simplex $\Delta_n$ is defined (for $n\in\mathbb N$) as a polytope with a set $X_n$ of $n+1$ affinely independent vertices, say $x_0,\dots,x_n$. Each element $x$ of $\Delta_n$ is expressed as a \emph{convex combination} or \emph{barycentric combination}
\begin{equation}\label{E:convcomb}
x = x_0 p_0 + \ldots + x_n p_n \ \mbox{with} \ p_i \in I \ \mbox{and} \  \sum_{i=0}^{n} p_i = 1
\end{equation}
of the vertices. The numbers $p_0,\dots,p_n$ are called the \emph{barycentric coordinates} of $x$. The affine independence of the vertices guarantees that each element $x$ of $\Delta_n$ has a unique representation as a convex combination of those vertices. An interior point of a face $F$ is a convex combination of the vertices belonging to the face in which each barycentric coordinate is positive. The face lattice of $\Delta_n$ is the power set of $X_n$.

It was shown by W. Neumann~\cite{N70} that the $n$-dimensional simplex $\Delta_n$ is the free barycentric algebra over the set $X_n = \{x_0,\dots,x_n\}$ of $n+1$ free generators.
This means that any mapping $f\colon X_n\to B$ from $X_n$ to any barycentric algebra $B$ has a unique extension to a barycentric homomorphism $\overline{f}\colon \Delta_n \to B$. The free barycentric algebra over the empty set is empty.

\begin{remark}
By~\cite[Lemma~5.8.1]{Modes},  elements of the free barycentric algebra on $n+1$ free generators $x_0, x_1, \dots, x_n$, may be expressed in the standard form $$(\dots((x_{0} x_{1} \underline{q}_1)x_{2} \underline{q}_2)\dots)x_{k}\underline{q}_k$$ for $q_i \in I$ and $k \leq n$, and by \cite[Lemma~1.6.4.1]{Modes}, the standard form may be reduced to the  convex combination of~\eqref{E:convcomb}. For a more direct connection, see \cite[\S2.5.1]{RSZ23}.
\end{remark}

Comparing~\eqref{E:convcomb} with~\eqref{E:affcomb}, we see that the vertices $x_0,\dots,x_n$ generate both $\Delta_n$ and the affine space isomorphic to $\mathbb{R}^n$, and that $\Delta_n$ is a subreduct (subalgebra of a reduct) of the affine space generated by $X_n$.
Note also that the affine space generated by the vertices of $\Delta_n$ is a free affine space over the same set of generators, and is isomorphic to $\mathbb{R}^n$. The combinations~\eqref{E:convcomb} and~\eqref{E:affcomb} do not depend on the enumeration of the generators, and for any point of the simplex the barycentric and affine coordinates coincide.

One model $S_n$ of a simplex $\Delta_n$ is especially useful, taking the free generators as the elements
\begin{equation}\label{E:freegens}
e_0 = (0,0,\dots,0, 0), e_1 = (1, 0,\dots,0,0), \dots, e_n = (0,0,\dots, 0, 1)
\end{equation}
of $\mathbb R^n$. These elements form an affinely independent set, and the simplex they generate is a free barycentric algebra over the set $\{e_0,\dots,e_n\}$ \cite{JK76}. The standard basis of the vector space $\mathbb{R}^n$ is then given by $\{e_1,\dots,e_n\}$.

\section{Barycentric coordinates in a polytope}\label{S:barcordpolytope}

Let $k\le n$ be natural numbers. Consider a polytope $P$ in $\mathbb{R}^k$ with vertex set $\mathbf v_0,\mathbf v_1,\dots,\mathbf v_n$. The polytope $P$ is a homomorphic image of the free barycentric algebra or simplex $\Delta_n$ on $n+1$ generators.  Thus, each element $\mathbf x$ of $P$ may be presented as a convex combination
\begin{align}\label{E:barcoord}
\mathbf{x}&=\mathbf v_0 p_0+\ldots+\mathbf v_n p_n
\\ \label{E:convcomb1}
\mbox{with}\quad
1&=p_0+\ldots+p_n
\end{align}
and $p_i \in I$. If $\mathbf x$ and $\mathbf v_i$ are given by Cartesian coordinates in the vector space $\mathbb{R}^k$, the barycentric coordinates $p_i$ may be calculated by solving the system \eqref{E:barcoord},\eqref{E:convcomb1} of $k+1$ linear equations in the $n+1$ nonegative unknowns $p_1,\dots,p_n$. If $P$ is not a simplex, then $k<n$, so the system is underdetermined, and the barycentric coordinates in \eqref{E:barcoord} are not uniquely specified.

The following problem appears in many applications of polytopes, for example in geometric modeling and computer graphics (see e.g. Floater~\cite{FloaterGBCA}).

\begin{problem}\label{P:problem}
Given the set $V$ of vertices $\mathbf v_i$ of a polytope $P$, find a uniform system to produce uniquely determined barycentric coordinates for any point $\mathbf x$ of $P$.
\end{problem}

Such a system must assign barycentric coordinates $p(\mathbf x,\mathbf v)$ to each point $\mathbf x$ of $P$ such that $\sum_{\mathbf v\in V} p(\mathbf x,\mathbf v) = 1$ and
\begin{equation}\label{E:applcomb}
\mathbf x = \sum_{\mathbf v\in V}p(\mathbf x,\mathbf v) \mathbf v,
\end{equation}
with some specific choice of $p(\mathbf x,\mathbf v) \in I$.
An extensive literature covers Problem~\ref{P:problem}, mostly written by mathematicians and engineers working on  applications of convex sets. (See e.g.~\cite{FloaterGBCA} and its references.) In \cite{RSZ23}, we discuss some of the available systems, in particular \emph{Gibbs coordinates} for barycentric algebras, based on entropy maximalization, and \emph{Wachspress coordinates} for convex polytopes. We also introduce \emph{chordal coordinates} for points in polygons based on chordal decompositions of polygons, together with a symmetrized version called \emph{cartographic coordinates}. In the present paper, we describe a system of sparse barycentric coordinates for polytopes. It is related to the chordal coordinates of polygons. However, it applies to general convex polytopes, and is more direct. On the other hand, it does not display the nice combinatorial properties of chordal coordinates. Our approach here is based on the fact that each polytope decomposes as a union of simplices, of the same dimension as the polytope in question.
We exploit the close relationship between presenting points of $\mbR^k$ as affine and as convex combinations of a given finite set of points.

Note that a \emph{decomposition} of a polytope $P$ of dimension $k$ is usually defined as a decomposition of $P$ into simplices of the same dimension which have pairwise disjoint interiors and have union equal to $P$. This is also known as a (generalized) \emph{triangulation}. There are numerous results concerning such triangulations. (See e.g.~\cite{LRS10}.)
In this paper we are interested in the following decomposition, which is possible for each polytope.

\begin{theorem}[Decomposition into simplices]~\cite[Thm.~3.1]{SRS09}\label{T:decomposition}
Let $P$ be a $k$-dimensional polytope with set $V$ of $n+1$ vertices. Fix $\mathbf{v} \in V$. Then $P$ is the union of simplices isomorphic to $\Delta_k$, which have pairwise disjoint interiors, each generated by a $(k+1)$-element subset of $V$ containing $\mathbf{v}$.
\end{theorem}

Denote the division of Theorem~\ref{T:decomposition} by $D_{v}^s$, and call it a \emph{$\mathbf v$-pointed decomposition}. The point $\mathbf v$ may be called a \emph{basic point} of $D_{v}^s$ or a \emph{basic vertex}. If we do not need to stress which vertex is basic, we just speak about \emph{pointed decompositions}.  Note that any two simplices of $D_{v}^s$ have a common wall that is a simplex containing $\mathbf v$.
The inductive proof of Theorem~\ref{T:decomposition} provided in \cite{SRS09} shows that, given pointed decompositions of facets (maximal walls) of $P$ not containing the vertex $\mathbf v$, each simplex of the decomposition $D_{v}^s$ of $P$ is generated by a simplex of the decomposition of a facet of $P$ and the vertex $\mathbf v$.

There is an intermediate step in proving Theorem~\ref{T:decomposition} which yields the existence of another interesting decomposition of the polytope $P$. Each facet $W$ of $P$, together with the basic point $\mathbf{v}$, forms a pyramid with basis $W$ and apex $\mathbf{v}$. The union of all these pyramids is $P$. Any two such pyramids have pairwise disjoint interiors, and their intersection is a pyramid with the apex $\mathbf{v}$. One obtains a decomposition $D_{v}^p$ of $P$ into the union of pyramids of the same dimension as $P$.

\begin{theorem}[Decomposition into pyramids]\label{T:decomposition2}
Let $P$ be a $k$-dimensional polytope with set $V$ of $n+1$ vertices. Fix $\mathbf{v} \in V$. Then $P$ is the union of those pyramids contained in $P$ which have a facet of $P$ as a basis, and the vertex $\mathbf{v}$ as the apex.
\end{theorem}

As in the case of decompositions of $P$ into a union of simplices, one also speaks about $\mathbf v$-pointed and pointed decompositions into pyramids. Note that a decomposition $D_{v}^s$ of $P$ may be obtained from a decomposition $D_{v}^p$ of $P$, using pointed decompositions of the pyramids into unions of simplices.

Now let us go back to the decomposition $D_{v}^s$ of $P$. Each point $\mathbf{x}$ of $P$ belongs to a simplex $S$ of $D_{v}^s$. Hence, $\mathbf{x}$ is the (uniquely defined) convex combination of some vertices of $S$. By taking the barycentric coordinates of the remaining vertices of $P$ to be $0$, one obtains the $\mathbf v$-\emph{pointed barycentric combination} of $\mathbf x$.

Note that the generators of $S$ also freely generate the affine space $\mathbb{R}^k$. Thus, one can also represent any point of $P$ (uniquely) as an affine combination of the vertices of $S$. However some of the corresponding affine coordinates may be negative.

To find $\mathbf v$-pointed barycentric coordinates of any point $\mathbf x$ of $P$, one needs a method of deciding to which simplex $S$ of $D_{v}^s$ the point $\mathbf x$ belongs. For example, if a line through $\mathbf v$ and $\mathbf x$ intersects a maximal wall of $S$ not containing $\mathbf v$, then the point $\mathbf x$ belongs to $S$, and can thus be presented as a (unique) convex combination of the vertices of $S$. However, we do not follow this idea in the current paper. In the following section, we exhibit other methods which work very nicely in the case of many polytopes, based on an appropriate numbering of the facets of $P$ and the simplices of its pointed decompositions.

\section{Pointed barycentric coordinates}\label{S:barcordDec}

We start with a well known and easy method of calculating the (unique) barycentric coordinates of points of a simplex, which will be frequently used later on.

\subsection{Volumetric coordinates in a simplex}\label{SS:volmtric}

The (unique) barycentric coordinates in a simplex may be given in terms of relative (hyper)volumes.

Consider a simplex $S$ of dimension $n$ that is spanned by the ordered set
$$
\set{\mathbf v_0,\dots,\mathbf v_{i-1},\mathbf v_i,\\ \mathbf v_{i+1},\dots,\mathbf v_n}
$$
in $n$-dimensional Euclidean space $\mathbb{R}^n$. Then each point $\mathbf x$ of $\mathbb{R}^n$ is uniquely described as an affine combination
\begin{equation}\label{E:affcsim}
\mathbf x = \sum_{i=0}^{n} \mathbf v_i p_i,
\end{equation}
where $p_0 = 1 - \sum_{i=1}^{n}p_{i}$.
If $\mathbf x$ and $\mathbf v_i$ are given by Cartesian coordinates of $\mathbb{R}^n$, the affine coordinates $p_1, \dots, p_n$ may be calculated by solving the equation \eqref{E:affcsim}. The unique solution is given by
\begin{equation}\label{E:volmtrpi}
p_i=\frac
{
\det[
\mathbf v_1-\mathbf v_0,
\dots,
\mathbf v_{i-1}-\mathbf v_0,
\mathbf x-\mathbf v_0,
\mathbf v_{i+1}-\mathbf v_0,
\dots,
\mathbf v_n-\mathbf v_0
]
}
{
\det[
\mathbf v_1-\mathbf v_0,
\dots,
\mathbf v_{i-1}-\mathbf v_0,
\mathbf v_i-\mathbf v_0,
\mathbf v_{i+1}-\mathbf v_0,
\dots,
\mathbf v_n-\mathbf v_0
]
}
,
\end{equation}
where $i = 1, \dots, n$. Now $\mathbf x$ belongs to the simplex $S$ precisely when all $p_i$ belong to the unit interval $I$. In this case, the $p_i$ are the (unique) barycentric coordinates of $\mathbf x$ with respect to $\mathbf v_i$ in the simplex.

If $\mathbf x \in S$ and $p_1, \dots, p_n \in I^{\circ}$, then $\mathbf x$ is an interior point of $S$. If all $p_i$ are in $I$ and at least one of $p_i$ is zero, then $\mathbf x$ belongs to the boundary of $S$, and when all but one of $p_i$ are zero, $\mathbf x$ is a vertex. If at least one of $p_i$ is negative, the point $\mathbf x$ does not belong to $S$.

Recall that the signed (hyper)volume of the simplex is equal to
\begin{equation}\label{E:hypvolv0}
(n!)^{-1}
\det[
\mathbf v_1-\mathbf v_0,
\dots,
\mathbf v_{i-1}-\mathbf v_0,
\mathbf v_i-\mathbf v_0,
\mathbf v_{i+1}-\mathbf v_0,
\dots,
\mathbf v_n-\mathbf v_0
].
\end{equation}
If $\mathbf x \in S$, a similar formula with $\mathbf x$ replacing $\mathbf v_i$ describes the (hyper)volume of the simplex spanned by $\mathbf v_0,\dots,\mathbf v_{i-1},\mathbf x,\mathbf v_{i+1},\dots,\mathbf v_n$. Then
the convex combination \eqref{E:affcsim} is said to give the \emph{volumetric} coordinates for $\mathbf x$ in terms of the vertices $\mathbf v_i$. (Cf. \cite[\S 4.1]{RSZ23}.)

\subsubsection{Areal coordinates}\label{SSS:Aov0v1v2}

The volumetric coordinates specialize, in the case of a triangle, to the classical \emph{areal coordinates}~\cite{Moebius},~\cite{Muggeridge}. The signed area of the triangle in the Euclidean plane $\mathbb R^2$ whose vertices are
$
\mathbf v_0=
\begin{bmatrix}
v_{00} &v_{01}
\end{bmatrix},
\mathbf v_1=
\begin{bmatrix}
v_{10} &v_{11}
\end{bmatrix}
\mbox{ and }
\mathbf v_2=
\begin{bmatrix}
v_{20} &v_{21}
\end{bmatrix}\,,
$
in counterclockwise order, is given by
\begin{equation}\label{E:Aov0v1v2}
A\left(\mathbf v_0,\mathbf v_1,\mathbf v_2\right)
=
\tfrac12\det[
\mathbf v_1-\mathbf v_0,
\mathbf v_2-\mathbf v_0
]
=
\frac12
\begin{vmatrix}
1 &v_{00} &v_{01}\\
1 &v_{10} &v_{11}\\
1 &v_{20} &v_{21}\\
\end{vmatrix}.
\end{equation}
Similar formulas hold for the signed areas
$A\left(\mathbf x,\mathbf v_1,\mathbf v_2\right)$, $A\left(\mathbf v_0,\mathbf x,\mathbf v_2\right)$ and $A\left(\mathbf v_0,\mathbf v_1,\mathbf x\right)$ spanned by an element $\mathbf x$ of the triangle and two of the vertices $\mathbf v_i$.

Note that for the area $A\left(\mathbf x,\mathbf y,\mathbf z\right)$ of a triangle spanned by counterclockwise ordered $\mathbf x < \mathbf y < \mathbf z$, one has
$$
A\left(\mathbf x,\mathbf y,\mathbf z\right) =
A\left(\mathbf y,\mathbf z,\mathbf x\right) =
A\left(\mathbf z,\mathbf x,\mathbf y\right).
$$

Moreover
\begin{equation}\label{E:arealcp1}
p_1=\frac
{
\det[
\mathbf x-\mathbf v_0,
\mathbf v_2-\mathbf v_0
]
}
{
\det[
\mathbf v_1-\mathbf v_0,
\mathbf v_2-\mathbf v_0
]
}
\,
= \frac
{
A\left(
\mathbf v_{0},\mathbf x,\mathbf v_{2}
\right)
}
{
A\left(
\mathbf v_{0},\mathbf v_{1},\mathbf v_{2}
\right)
},
\end{equation}
\begin{equation}\label{E:arealcp2}
p_2=\frac
{
\det[
\mathbf v_1-\mathbf v_0,
\mathbf x-\mathbf v_2
]
}
{
\det[
\mathbf v_1-\mathbf v_0,
\mathbf v_2-\mathbf v_0
]
}
\,
= \frac
{
A\left(
\mathbf v_{0},\mathbf v_{1},\mathbf x
\right)
}
{
A\left(
\mathbf v_{0},\mathbf v_{1},\mathbf v_{2}
\right)
},
\end{equation}
and
\begin{equation}\label{E:arealcp0}
p_0 = 1 - p_1 - p_2 =
 \frac
{
A\left(
\mathbf x,\mathbf v_{1},\mathbf v_{2}
\right)
}
{
A\left(
\mathbf v_{0},\mathbf v_{1},\mathbf v_{2}
\right)
},
\end{equation}

\subsection{Pointed barycentric coordinates in polygons}\label{S:barcordpolgs}

Now suppose that $\Pi$ is a polygon spanned by a counterclockwise ordered $n$-tuple of vertices $\mathbf v_{1} < \dots < \mathbf v_{n}$. Consider the decomposition $\mathcal{D}_1 = \mathcal{D}_{v_1}$ of $\Pi$ into the union of the triangles $\tau_{1i (i+1)}$ given by Theorem~\ref{T:decomposition}, where $\tau_{1i (i+1)}$ denotes the triangle spanned by $\mathbf v_{1} < \mathbf v_{i} < \mathbf v_{i+1}$ with $i = 2, \dots, n-1$. The triangles form the sequence
\begin{equation}\label{E:sequence}
\tau_{123}, \tau_{134}, \dots, \tau_{1 (n-1) n}.
\end{equation}
A point $\mathbf{x}$ of $\Pi$ which is not itself a vertex belongs either to precisely one of the triangles, or to a common side of two adjacent triangles. To decide to which triangle of \eqref{E:sequence} the point $\mathbf{x}$ belongs, one calculates the affine coordinates $p_1, p_i, p_{i+1}$ of $\mathbf{x}$ in each of the triangles $\tau_{1i (i+1)}$, starting from the first one, then in the order of~\eqref{E:sequence}. (In fact, it suffices to consider the signs of the numerators.) The first $\tau_{1i (i+1)}$ with $0 \leq p_1, p_i, p_{i+1} < 1$ and at least one of $p_1, p_i, p_{i+1}$ positive is the triangle containing $\mathbf{x}$. The $\mathbf v_1$-pointed barycentric coordinates $p_1, p_i, p_{i+1}$ of $\mathbf{x}$ are given by \eqref{E:arealcp0}, \eqref{E:arealcp1}, and \eqref{E:arealcp2}, with $1$ replacing $0$ and $i$ replacing $1$, and all other coordinates equal to $0$.

\begin{example}\label{Ex: 6-gon}

The picture
$$
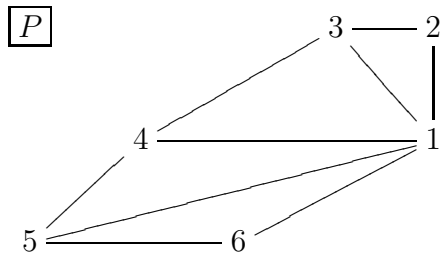
\xymatrix{
\boxed{P}
&
&
&
3
\ar@{-}[r]
\ar@{-}[dr]
\ar@{-}[dll]
&
2
\\
&
4
&
&
&
1
\ar@{-}[lll]
\ar@{-}[dllll]
\ar@{-}[u]
\ar@{-}[dll]
\\
5
\ar@{-}[ur]
\ar@{-}[rr]
&
&
6
}
$$
illustrates a $1$-pointed decomposition of a hexagon $P$.
\end{example}

The pointed decompositions of a polygon $\Pi$ are special chordal decompositions. In \cite[\S~7]{RSZ23}, chordal decompositions of a polygon are defined as follows.
The \emph{skeleton} of the polygon $\Pi$ is the cyclic graph
$C_n$ constituted by the vertices and undirected edges of the polygon. In the cyclic graph $C_n$, a \emph{chord} is an edge connecting vertices which are not adjacent in $C_n$.
A \emph{chordal decomposition} of the polygon $\Pi$ with ordered vertex set $V = \{\mathbf v_1<\mathbf v_2<\dots<\mathbf v_n\}$ is a system of $n-3$ non-crossing chords of $C_n$ that decompose $\Pi$ as a union of $n-2$ simplices (triangles) whose vertices are vertices of $\Pi$.
Given a chordal decomposition, one obtains others by the action of the dihedral group $D_n$.
If $n < 6$, the pointed decompositions of $\Pi$ comprise all of its chordal decompositions. However, for $n > 5$, there are also non-pointed chordal decompositions of $\Pi$ into a union of triangles \cite[\S 7]{RSZ23}.

While pointed barycentric coordinates of elements of $\Pi$ are sparse, in that they involve minimal numbers of vertices, more symmetric coordinatizations are given by \emph{cartographic coordinates} obtained as an average of $\mathbf v$-pointed barycentric coordinates over all vertices $\mathbf v$ of $\Pi$ \cite[\S~7]{RSZ23}.

\subsection{Pointed barycentric coordinates in polytopes}

Now consider a polytope $P$ of any higher dimension, with its given sets of vertices and facets. To find pointed barycentric coordinates of a point $\mathbf x$ of $P$, one needs first to decompose $P$ into a union of pyramids, as described in Theorem~\ref{T:decomposition2}, and then to decompose each pyramid into a union of simplices, as described in Theorem~\ref{T:decomposition}. To find the simplex containing a given point $\mathbf x$, it is convenient to keep track of the facets and the simplices of the decomposition, along with their orientations.

While for polygons the simplices of pointed decompositions are readily numbered in consecutive order, the situation becomes much more complicated in the case of polytopes of higher dimension. A naive method may be based on first numbering facets of a polytope, and then numbering the vertices inside of each facet following the numbering of the facets. However, this should be done in such a way that the labels given to vertices will not be repeated, and such that adjacent facets will be numbered consecutively. With the facets and vertices numbered, one may linearly order and number the simplices of the decomposition. Then, one subsequently calculates the affine coordinates of a given point $\mathbf x$ in the numbered simplices, until arriving at a simplex $S$ in which the affine coordinates of $\mathbf x$ are barycentric \textemdash i.e., all are non-negative. The (volumetric) coordinates of $\mathbf x$ in $S$, together with the coordinates $0$ for vertices not contained in $S$, provide the pointed barycentric coordinates of $\mathbf x$ in $P$.

The last part of this section is devoted to the problem of numbering and orienting the facets, vertices, pyramids and simplices in a pointed decomposition of a polytope, especially in the case of $3$-dimensional polytopes.

\subsubsection{Pointed barycentric coordinates in $3$-dimensional polytopes}

First note that each (non-empty) system $\mathcal{F}$ of facets of any $k$-dimensional polytope $P$ may be numbered using the method called a \emph{semi-shelling}, as described in \cite[Ex~15.1]{AB83}. A \emph{semi-shelling} of $P$ is a numbering $F_1, \dots, F_r$ of the facets of $\mathcal{F}$ such that for $i = 2, \dots, r$,
\begin{equation}\label{E:shell}
F_i \cap \bigcup_{j=1}^{i-1} F_j
\end{equation}
is a non-empty union of $(k-2)$-dimensional faces of $F_i$. Note that any facet of $P$ can initially be taken as $F_1$. A semi-shelling is a \emph{shelling} if, additionally, the set of~\eqref{E:shell} is homeomorphic to a $(k-2)$-ball or $(k-2)$-sphere \cite{BayerLee}. Note that every polytope is shellable \cite[Th.~2.1]{BayerLee}.

Let $\Pi$ be a $3$-dimensional polytope with a set $\mathbf V$ of vertices. Our aim is to decompose $\Pi$ into a union of simplices, and then to number these simplices. First recall that all facets of $\Pi$ are polygons, and we may assume that the vertices of each of them are counterclockwise ordered when started from any one of them.

Our method of decomposing $\Pi$ into a union of simplices and then numbering them is described in the following steps.

\begin{enumerate}

\item Fix a vertex $\mathbf v$ of $\mathbf V$ as the basic point and denote it by $0$. The remaining vertices will be denoted by positive integers.

\item Using the shelling method, first number the facets of $\Pi$ containing $0$, and denote them by $F_1, \dots, F_r$.

\item Continue to number the remaining facets of $\Pi$ (not containing $\mathbf v$), using the shelling method again, and denote them consecutively by $F_{r+1}, \dots, F_{l}$.

\item Each facet $F_j$, for $j = 1, \dots, l$, is the polygon spanned by the counterclockwise ordered $n_j$-tuple
\begin{equation}
1_j < 2_j < \dots < n_j
\end{equation}
of its vertices. Any vertex of $F_j$ may be chosen as the first one.

\item The vertices of $\Pi$ are numbered as follows. At first, the vertices $1_1, \dots, n_1$ of $F_1$ are numbered as $1, \dots, n_1$. For $j > 1$, if a vertex $s_j$ of $F_j$ appears as $p_i$ in $F_i$ for some $i < j$, then $s_j = p_i$. New numbers are attached consecutively only to vertices of $F_j$ which are not members of $F_i$ for some $i < j$.
Consequently, the vertices of $\Pi$ may be numbered and ordered as follows:
\begin{equation}\label{E:numbvert}
0 < 1 \dots < n_1 < 1'_{r+1} < \dots n'_{r+1} < \dots < 1'_l < \dots < n'_{l},
\end{equation}
where $1'_j < \dots < n'_j$ is the subsequence of $1_j < \dots < n_j$ obtained by removing  all numbers which appeared already earlier in~\eqref{E:numbvert}. Note that it is possible that all members of $1_j < \dots < n_j$ will disappear in this way. The vertices of
~\eqref{E:numbvert} may now be renumbered with consecutive natural numbers.

\item With facets and vertices already numbered, proceed to numbering the simplices of the pointed decompositions of the facets. Decompose each facet $F_j$ for $j = r+1, \dots, l$, taking $1_j$ as a basic point, into the union of triangles $S_{j_i} = \tau_{1_{j} i_{j} {(i+1)}_{j}}$ with $i = 2, \dots, n$, and form the sequence $S_{1_i}, \dots, S_{n_j}$. All the simplices obtained in this way from the facets $F_j$ form the sequence
\begin{equation}\label{E:numbsimpl}
S_{1_1}, \dots, S_{1_{n_j}}, \dots, S_{l_1}, \dots, S_{n_l},
\end{equation}
which is now numbered by consecutive natural numbers.

\item Each simplex $\mathcal{S}_{j_i}$ of the $0$-pointed decomposition $D_{0}^s$ of $\Pi$ is then obtained as the cone generated by the basic vertex $0$ and the simplex $S_{j_i}$ for $j = r+1, \dots, l$ and $i = 1, \dots, n$. The numbering of the simplices $\mathcal{S}_{j_i}$ coincides with that given to the simplices $S_{j_i}$.
\end{enumerate}

\begin{example}
Consider the hexagon $P$ of Example~\ref{Ex: 6-gon}, and the pyramid $P_6$ with the hexagon $P$ as its basis, and a new vertex $0$ at the apex:
$$
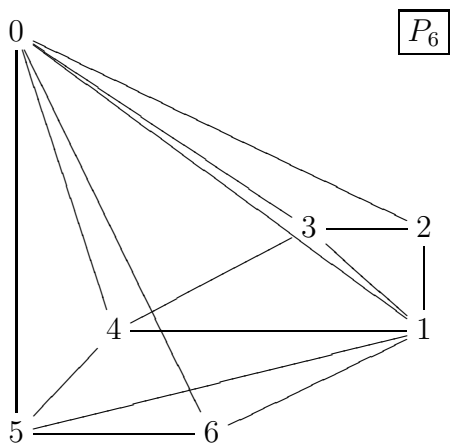
\xymatrix{
0
\ar@{-}[dddd]
\ar@{-}[dddr]
\ar@{-}[ddddrr]
\ar@{-}[dddrrrr]
\ar@{-}[ddrrr]
\ar@{-}[ddrrrr]
&&&&
\boxed{P_6}
\\
\\
&
&
&
3
\ar@{-}[r]
\ar@{-}[dr]
\ar@{-}[dll]
&
2
\\
&
4
&
&
&
1
\ar@{-}[lll]
\ar@{-}[dllll]
\ar@{-}[u]
\ar@{-}[dll]
\\
5
\ar@{-}[ur]
\ar@{-}[rr]
&
&
6
}
$$
Then $P$ is the only facet of $P_6$ not containing $0$. Now take the pointed decomposition of $P$ given in Example~\ref{Ex: 6-gon}. The facets of $P$  provide the four triangles $S_1 = \tau_{123}, S_2 = \tau_{134}, S_3 = \tau_{145}, S_4 = \tau_{156}$ of the decomposition of $P$. They determine the four simplices $\mathcal{S}_1, \mathcal{S}_2, \mathcal{S}_3, \mathcal{S}_4$ of the $0$-pointed decomposition of $P_6$, obtained as  the cones generated by the simplices $S_1, S_2,  S_3, S_4$ and the vertex $0$.
\end{example}

\begin{figure}[hbt]
\centering
\caption{The lopped cube of Example~\ref{X:LoppdCub}.}
\label{F:LoppdCub}
\scalebox{0.9}
{
\begin{tikzpicture}
    \coordinate (Origin)   at (0,0);
    \coordinate (XAxisMin) at (-5,0);
    \coordinate (XAxisMax) at (5,0);
    \coordinate (YAxisMin) at (0,-5);
    \coordinate (YAxisMax) at (0,5);
    \clip (-5.5,-5.5) rectangle (5.5cm,5cm);
\node[label=left:{$7$}, draw, circle, inner sep=2pt, fill] at (-1,4) {};
\node[label=right:{$5$}, draw, circle, inner sep=2pt, fill] at (1,4) {};
\node[label=above:{$F_3$}] at (0,3.2) {};
\node[label=above:{$F_2$}] at (2,1.2) {};
\node[label=above:{$F_4$}] at (-2,1.2) {};
\node[label=above:{$F_6$}] at (2,-2) {};
\node[label=above:{$F_7\colon 45789$}] at (3.6,-4) {};
\node[label=above:{$F_5$}] at (-2,-2) {};
\node[label=above:{$F_1$}] at (0,-.4) {};
\node[label=left:{$6$}, draw, circle, inner sep=2pt, fill] at (0,3) {};
\node[label=right:{$0$}, draw, circle, inner sep=2pt, fill] at (0,1.5) {};
\node[label=right:{$2$}, draw, circle, inner sep=2pt, fill] at (0,-1.5) {};
\node[label=left:{$8$}, draw, circle, inner sep=2pt, fill] at (-5,0) {};
\node[label=right:{$1$}, draw, circle, inner sep=2pt, fill] at (-1.5,0) {};
\node[label=left:{$3$}, draw, circle, inner sep=2pt, fill] at (1.5,0) {};
\node[label=right:{$4$}, draw, circle, inner sep=2pt, fill] at (5,0) {};
\node[label=below:{$9$}, draw, circle, inner sep=2pt, fill] at (0,-5) {};
\draw (-1,4) -- (1,4);
\draw (-1,4) -- (0,3);
\draw (0,1.5) -- (0,3);
\draw (1,4) -- (0,3);
\draw (-1,4) -- (-5,0);
\draw (1,4) -- (5,0);
\draw (-5,0) -- (-1.5,0);
\draw (-5,0) -- (0,-5);
\draw (5,0) -- (1.5,0);
\draw (5,0) -- (0,-5);
\draw (0,1.5) -- (-1.5,0);
\draw (0,-1.5) -- (0,-5);
\draw (-1.5,0) -- (0,-1.5);
\draw (-1.5,0) -- (0,1.5);
\draw (1.5,0) -- (0,-1.5);
\draw (1.5,0) -- (0,1.5);
\end{tikzpicture}
}
\end{figure}
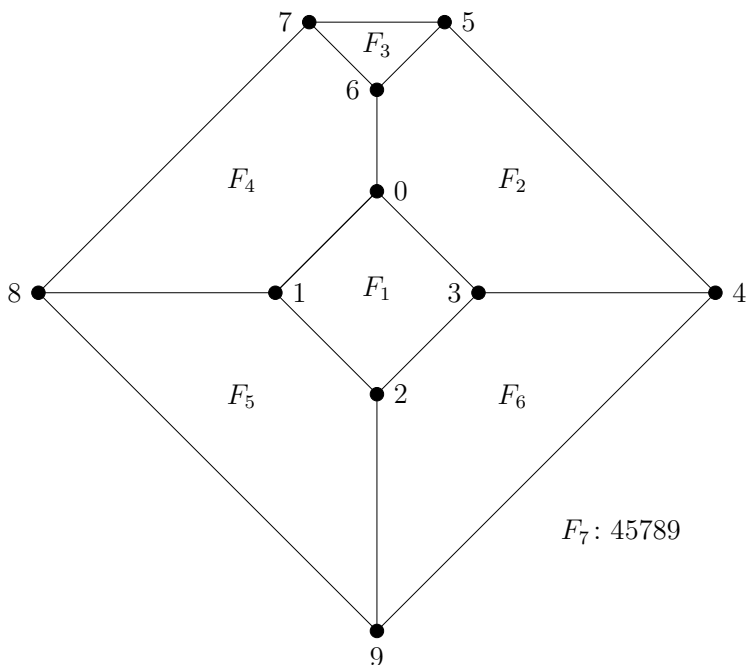

\begin{example}\label{X:LoppdCub}
The skeleton of a $3$-dimensional polytope $\Pi$ is a planar graph (compare \cite[Ex~15.2]{AB83}), so that it is possible to sketch such polytopes on the plane. For example, consider the projection of the lopped cube $C$  on the plane as in Figure~\ref{F:LoppdCub}, and
its seven facets.

The numbering of the facets is one of the possible numberings obtained by using the shelling method. The numbering of the vertices follows the numbering of the facets. Then $F_3, F_5, F_6, F_7$ are the facets not containg the vertex $0$.  To obtain a numbering of the simplices of the $0$-pointed decomposition of $C$, one decomposes each of these facets into a union of triangles. The simplices of the $0$-pointed decomposition of $C$ are the simplices obtained as the cones generated by those triangles and the vertex $0$.
\end{example}

To find barycentric coordinates of a given point $\mathbf x$ of $\Pi$, one calculates the affine coordinates of $\mathbf x$ in all the simplices of the decomposition $D_{0}^s$ of $\Pi$, in the order  described above, until one finds a simplex $\mathcal{S}$ for which the affine coordinates of $\mathbf x$ are barycentric, similarly as was done in the case of polygons. If all these coordinates are positive, then $\mathbf x$ is an interior point of $\mathcal{S}$. If some of them are $0$, then $\mathbf x$ belongs to the boundary of $\mathcal{S}$. The barycentric coordinates of the point $\mathbf x$ are given by the volumetric coordinates of $\mathbf x$ in $S$, with all other coordinates equal to $0$.

\subsubsection{Pointed barycentric coordinates in $n$-dimensional polytopes}

The method used to find barycentric coordinates of points in a $3$-dimensional polytope extends recursively to polytopes of higher dimensions.

\end{document}